\documentclass[12pt]{amsart}
\usepackage[margin=3cm]{geometry}
\RequirePackage{amsmath, amssymb, amsthm, amsfonts, amscd}
\RequirePackage{indentfirst}
\RequirePackage[pdftex]{graphicx}
\RequirePackage{hyperref}
\RequirePackage[dvipsnames]{xcolor}
\usepackage{algorithm}
\usepackage{algorithmic}

\numberwithin{equation}{section}

\newtheorem{theorem}{Theorem}[section]
\newtheorem{corollary}[theorem]{Corollary} 
\numberwithin{claim}{section}

\newtheorem{remark}[theorem]{Remark}  

\theoremstyle{definition}
\numberwithin{definition}{section}
\newtheorem*{definition*}{Definition}

\numberwithin{example}{section}

\newcommand{\Sp}{ \hspace{0.05cm} }

\usepackage{amsopn}

\begin{document}

\title[Carleman-Fourier linearization]{Carleman-Fourier linearization of nonlinear real dynamical systems }

\author{ Nader Motee and Qiyu Sun}
\address{Motee: Mechanical Engineering and Mechanics, Lehigh University,  Bethlehem, PA 18015}
\email{nam211@lehigh.edu}
\address{Sun: Department of Mathematics, University of Central Florida, Orlando, Florida 32816}
\email{qiyu.sun@ucf.edu}

\date{}

\begin{abstract}
This paper presents Carleman-Fourier linearization
 for analyzing nonlinear real dynamical systems with periodic vector fields. 
Using Fourier basis functions, this novel framework transforms such dynamical systems into equivalent infinite-dimensional linear  dynamical systems. In this paper, we establish the exponential convergence of the primary block in the finite-section approximation of this linearized 
system to the state vector of the original nonlinear system. To showcase the efficacy of our approach, we apply it to the Kuramoto model, a prominent model for coupled oscillators. The results demonstrate promising accuracy in approximating the original system's behavior.
\end{abstract}


\maketitle

\section{Introduction}

Nonlinear real dynamical systems with (quasi-)periodic vector fields are prevalent across various scientific disciplines, including biological systems, physical systems, and engineering systems. The inherent nonlinearity 
of these systems presents a significant challenge in developing a comprehensive mathematical framework for their analysis. One mainstream approach to studying nonlinear dynamical systems is to transform them into linear counterparts, leveraging well-established techniques for linear systems. Carleman linearization, introduced in 1932, is a prominent method to transform finite-dimensional nonlinear systems into infinite-dimensional linear systems.
Similar to the Maclaurin expansion for analytic functions, Carleman linearization is particularly effective for systems whose dynamics can be well approximated by low-degree polynomials, and the
resulting infinite-dimensional linear system  can be   exponentially approached by the finite-section method 
\cite{ Amini2022, brockett2014early, Brunton2016, Foret2017, forets2021reachability, Korda2018, KordaMezic2018, Kowalski1991, minisini2007carleman, steeb1980non, WangJungersOng2023}.

Carleman linearization,  achieved through state variable monomials,
has seen a recent resurgence due to advances in theoretical understanding, improved numerical and algorithmic techniques, and increased availability of extensive datasets \cite{Akiba2023, Amini2022, Brunton2016,  Liu2021}. Despite its theoretical advantages,  it remains a challenge to achieve well representations for those nonlinear systems poorly approximated by polynomials. In this paper, we present a novel approach that leverages the Fourier representation of periodic vector fields in conjunction with traditional Carleman linearization techniques. By exploiting the inherent structure of dynamical systems with periodic vector fields, our method  captures complicated nonlinear behaviors more efficiently,  enhances the parsimony and interpretability of the resulting embedding,
and surpasses the accuracy of existing linearization techniques like Carleman  linearization using monomials.

\smallskip

Let us start from considering an illustrative nonlinear dynamical system on the real  line \(\mathbb{R}\),  
\begin{equation}\label{dynamicsystem.def}
\dot{x} = g(x) \ \ {\rm with}\  \ x(0) = x_0,
\end{equation}
 governed by
 a real-valued periodic function  \(g(x)\) that can be expressed as a Fourier series:
\begin{equation}\label{Fourierexpansion.def}
g(x) = \sum_{n \in \mathbb{Z}} g_n e^{i n x},
\end{equation}
where \(\mathbb{Z}\) denotes the set of all integers and $g_n, n\in {\mathbb Z}$, are the Fourier coefficients. These coefficients decay exponentially; that is, there exist constants \(D > 0\) and \(0 < r < 1\) such that
\begin{equation}\label{fouriercoefficient.decay}
|g_n| \le D r^{|n|}\ \ \textrm{for  all}\  n \in \mathbb{Z}.
\end{equation}
This implies that $g(x)$ is analytic in a neighborhood of the origin and can be expressed as a Maclaurin series:
\begin{equation}\label{g.Taylor}
g(x) = \sum_{m \in \mathbb{Z}_+} c_m x^m \quad \text{for} \quad x \in \left( -\ln\frac{1}{r}, \ln\frac{1}{r} \right),
\end{equation}
 where $\mathbb{Z}_+$ denotes the set of all nonnegative integers.
  By applying the traditional Carleman linearization  \eqref{Carleman.eq}—which involves substituting the Maclaurin expansion of $g(x)$
 from  \eqref{g.Taylor} into the dynamical system in  \eqref{dynamicsystem.def}—it has been shown in \cite{Amini2022} that under the conditions
\begin{equation}\label{Carleman.requirement}
g(0) = 0 \quad \text{and} \quad |x_0| < e^{-1} \ln\left( \dfrac{1}{r} \right),
\end{equation}
the principal component of the finite-section approach for the infinite-dimensional linear system in ~\eqref{Carleman.eq} converges exponentially to the state $x(t)$ of the original nonlinear system over a short time period, as demonstrated in ~\eqref{Carleman.exponential}.

\smallskip

The requirements specified in \eqref{Carleman.requirement} are not universally satisfied. These conditions imply that the origin must serve as an equilibrium point for the dynamical system in \eqref{dynamicsystem.def}, and that the initial state must be close to this equilibrium. This limitation restricts the broader applicability of the Carleman linearization. To alleviate these constraints, we redefine the state variables using complex exponentials of extended state variables $\tilde{\mathbf{x}} = [x, -x]^{\mathsf{T}}$, which includes both the state variable and its negative, where $\mathsf{T}$ denotes the vector transpose. This approach allows us to propose a novel Carleman-Fourier linearization of the nonlinear dynamical system in \eqref{dynamicsystem.def}, where the state matrix takes the form of a block-upper triangular matrix; see \eqref{CarlemanFourier.def}. This block-upper triangular structure is essential for deriving explicit error bounds and proving convergence results. Furthermore, we demonstrate that the finite-section approach to the proposed Carleman-Fourier linearization  \eqref{CarlemanFourier.def} achieves exponential convergence over a specified time range as the model size increases, {\bf without} requiring the conditions in \eqref{Carleman.requirement} for the equilibrium and initial state of the dynamical system in \eqref{dynamicsystem.def}; see Theorem \ref{main.thm} and Corollary \ref{main.cor1}.

Dynamical systems with quasi-periodic vector fields, characterized by multiple incommensurate fundamental frequencies, exhibit remarkable complexity in their internal dynamics. The interactions among these frequencies generate a wide range of behaviors, including intricate phase space structures, frequency locking phenomena, and even chaotic dynamics \cite{Broer1996, gentile2010, glazier1998, liao2020, suzuki2016}. In this paper, we extend the Carleman-Fourier linearization method   from the illustrative real dynamical systems described in \eqref{dynamicsystem.def} to more general systems in \eqref{dynamicsystem.multiple}, where the vector field is a quasi-periodic vector-valued function with multiple fundamental frequencies.  We establish exponential convergence results for its finite-section approach, which could pave the way for efficient reachability analysis and novel computational algorithms for quasi-periodic nonlinear systems.

The Carleman-Fourier linearization offers several key advantages over the classical Carleman linearization, particularly in handling periodic vector fields more accurately over larger neighborhoods. This capability is crucial for analyzing system behavior beyond equilibrium points, enabling reliable predictions under significant perturbations and varying conditions. Additionally, the resulting linear system improves predictability over longer time intervals, which is essential for studying the long-term behavior of these systems. This improvement arises from the efficiency of Fourier basis functions in capturing both periodic and nonlinear dynamics more effectively than polynomial bases.

Carleman–Fourier linearization also offers significant advantages for quasi-periodic systems by generating linear systems with parsimonious matrices, owing to the inherent sparsity of their Fourier representations. This sparsity reduces computational costs and simplifies optimization processes, thereby enhancing the scalability of reachability analysis and learning algorithms that utilize such linear approximations \cite{Brunton2016}. Since quasi-periodic functions often possess sparse Fourier representations involving only a few fundamental frequencies, the complexity of analyzing these systems is substantially diminished. This reduction opens new opportunities to effectively learn and model this class of systems. Advances in compressed sensing and sparse regression further facilitate the identification of non-zero Fourier coefficients using data-driven algorithms, eliminating the need for combinatorially expensive searches \cite{Brunton2016,  Kramer2024, Yu2024}. By focusing on the most relevant components of the system, these techniques reduce unnecessary complexity and ensure practical and efficient modeling. Ultimately, Carleman–Fourier linearization achieves a balance between model simplicity and computational efficiency, resulting in models that are interpretable, scalable, and suitable for real-world applications.

\section{Classical Carleman Linearization}
\label{Carleman.section}

In the Maclaurin series \eqref{g.Taylor}, the coefficients \( c_m,  m \in \mathbb{Z}_+ \),  are expressed in terms of the Fourier coefficients \( g_n \) as:
\[
c_m = \frac{1}{m!} \sum_{n \in \mathbb{Z}} (i n)^m g_n.
\]
The above Maclaurin  expansion applies  since
\begin{equation*} |c_0|=|g(0)|\le  \sum_{n\in {\mathbb Z}} |g_n|\le D+ \frac{2D r }{1-r}= D \left(\frac{1+r}{1-r}\right)
\end{equation*}
and for $m\ge 1$
\begin{equation*}
|c_m| \le     \frac{2 D}{m!}   \sum_{n=1}^\infty n^m r^n
\le \frac{2D}{r} \int_1^\infty \frac{t^m r^t}{m!} dt
\le   \frac{2D}{r} \left(\ln \left( \frac{1}{r}\right)\right)^{-m-1}.\end{equation*}
Expanding \(g(x)\) as a Maclaurin series yields the classical Carleman linearization  of the nonlinear dynamical system described by Eq.~(\ref{dynamicsystem.def}) \cite{Amini2022, Foret2017, Kowalski1991, Liu2021}. This linearization can be expressed as
\begin{equation}\label{Carleman.eq}
\dot{\mathbf{x}} = \mathbf{A}\mathbf{x} + \mathbf{a},
\end{equation}
where \(\mathbf{x} = [x, x^{2}, \ldots, x^{N}, \ldots]^{\mathsf{T}}\), \(\mathbf{a} = [c_{0}, 0, \ldots, 0, \ldots]^{\mathsf{T}}\), and the state matrix \(\mathbf{A} = [n c_{n'-n+1}]_{n,n'=1}^{\infty}\) is given by
\begin{equation}\label{Carlemanmatrix.def-A}
\mathbf{A} =
\begin{bmatrix}
 c_{1} &  c_{2} & \cdots & c_{n-1} & c_{n} & \cdots \\[6pt]
 2c_{0} & 2c_{1} & \cdots & 2c_{n-2} & 2c_{n-1} & \cdots \\[6pt]
   0     & \ddots & \ddots & \vdots & \vdots & \ddots \\[6pt]
    0    &    0    & \ddots & \ddots & \vdots & \ddots \\[6pt]
     \vdots   &        &    \ddots    & n c_{0} & n c_{1} & \ddots \\[6pt]
        &        &        &         & \ddots  & \ddots
\end{bmatrix}.
\end{equation}
The classical Carleman linearization is based on decoupling the monomials $x^n, n\ge 1$, and it can be viewed as an {\bf analogue} of the Maclaurin expansion for analytic functions in the setting of dynamical systems. Since the state matrix $ {\bf A} $ is an unbounded operator on the space of square-summable sequences, we consider the following finite-dimensional dynamical systems of size $ N \ge 1 $,
\begin{equation}\label{finitesection.Carleman}
 \dot {\bf x}_N=   \begin{bmatrix}
 c_1 &  c_2 & \dots & c_{n-1} & c_n \\
     2c_0 & 2c_1 & \cdots  & 2c_{n-2} &  2 c_{n-1}\\
     &\ddots & \ddots & \vdots & \vdots\\
        & & \ddots & \ddots & \vdots \\
      & &  & n c_0&  n c_1 \\
\end{bmatrix}{\bf x}_N + \begin{bmatrix}
    c_0 \\
    0 \\
    \vdots\\ \vdots \\ 0\\
\end{bmatrix},
\end{equation}
with the initial condition $ {\bf x}_N(0)= [x_0, x_0^2, \ldots, x_0^N]^{\mathsf T}$. These systems are frequently employed to approximate the infinite-dimensional linear dynamical system in \eqref{Carleman.eq}, where $ {\bf x}_N = [x_{1, N}, \ldots, x_{N, N}]^{\mathsf{T}} $, and the state matrix is the leading $ N \times N $ principal submatrix of $ {\bf A} $. Under the assumption in  \eqref{Carleman.requirement}, and using the upper-triangular structure of $ {\bf A} $, we follow the argument of  \cite{Amini2022} to show that the first component $ x_{1, N} $ of the state vector $ {\bf x}_N $ in \eqref{finitesection.Carleman} provides an exponential approximation to the state $ x $ of the original nonlinear dynamical system \eqref{dynamicsystem.def} over a certain time interval. Specifically, there exist $ T^* > 0 $ and absolute constants $ C > 0 $ and $ r_1 \in (0,1) $, independent of $N$, such that
\begin{equation}\label{Carleman.exponential}
|x_{1, N}(t)-x(t)| \le C r_1^N
\end{equation}
hold for all $ 0 \le t \le T^* $ and $ N \ge 1 $.

The conditions in \eqref{Carleman.requirement} that ensure exponential convergence described in \eqref{Carleman.exponential} are not universally met, especially regarding the equilibrium point criterion $ g(0)=c_0=0 $. In principle, if this equilibrium condition is relaxed to allow for $ g(0) $ to be sufficiently close to the origin, then the finite-section approximation to the Carleman linearization should still exhibit exponential convergence.
Numerical evidence on simulations of the Kuramoto model supports this reasoning; see the top right plot of Fig. \ref{fig:kuramoto}, and it indicate that when the initial phase is near zero, the finite-section approach attains exponential convergence, whereas the convergence fails if the initial phase does not remain small.

Existing proofs of exponential convergence, such as those in \cite{Amini2022, Foret2017, Liu2021}, rely heavily on the upper-triangular structure of the state matrix $ {\bf A} $ in \eqref{Carlemanmatrix.def-A}. When the equilibrium condition in \eqref{Carleman.requirement} is not met, this structural property is lost and the approach in \cite{Amini2022, Foret2017, Liu2021} cannot be applied directly.
 This limitation restricts the applicability of standard Carleman linearization to nonlinear systems with periodic right-hand sides. By contrast, the Carleman–Fourier linearization introduced in the next section addresses this difficulty. It ensures that the associated finite-section approximations achieve exponential convergence, and the convergence rate can be selected to be {\bf independent} of both $ g(0)=c_0 $ and the initial condition $ x_0 $, as established in Theorem \ref{main.thm} and Corollary \ref{main.cor1}.

\section{Carleman-Fourier linearization and its finite-section approximation}
\label{CarlemanFourier.section}

The selection of an appropriate basis for linearization is critical to ensuring that the resulting state matrix is both analytically tractable and computationally efficient. While the polynomial basis \(\{x^n\}_{n \geq 1}\) is traditionally employed in the Carleman linearization, its application to dynamical systems involving the periodic function \(g(x)\) in \eqref{dynamicsystem.def} proves mathematically incongruous. Although the resulting state matrix \({\bf A}\) derived using the polynomial basis is indeed upper triangular—a property that facilitates
the exponential convergence of the finite-section method in \eqref{Carleman.exponential}—its upper triangular part is not sparse. This lack of sparsity imposes significant computational burdens and limits scalability to higher-dimensional problems. Our overarching aim, therefore, is to construct a linearization that not only preserves the upper triangular structure but also promotes sparsity within that structure.

To fulfill this objective, adopting the Fourier basis \(\{e^{inx}\}_{n \in \mathbb{Z}_0}\), as formulated in \eqref{CarlemanFourier.eq0}, better exploits the periodicity of \(g(x)\). Specifically, let \(x\) be as given in \eqref{dynamicsystem.def}, and define
$
{\bf w} = [e^{inx}]_{n \in \mathbb{Z}_0},
$
where \(\mathbb{Z}_0\) denotes the set of all nonzero integers. It follows that
\begin{equation}\label{CarlemanFourier.eq0}
\dot{\bf w} = {\bf H}\,{\bf w} + {\bf h},
\end{equation}
where
the state matrix \({\bf H} = [in g_{n'-n}]_{n, n' \in \mathbb{Z}_0}\) and
nonhomogenous term
\({\bf h} = [\,in\,g_{-n}\,]_{n \in \mathbb{Z}_0}\)
are given by
\begin{equation}\label{Carlemanmatrix.def}
{\bf H} =
\begin{bmatrix}
\ddots & \vdots & \vdots & \vdots & \vdots & \ddots\\
\cdots & -2ig_0 & -2ig_{1} & -2ig_{3} & -2ig_{4} & \cdots\\
\cdots & -ig_{-1} & -ig_0 & -ig_{2} & -ig_{3} & \cdots\\
\cdots & ig_{-3} & ig_{-2} & ig_0 & ig_{1} & \cdots\\
\cdots & 2ig_{-4} & 2ig_{-3} & 2ig_{-1} & 2ig_0 & \cdots\\
\ddots & \vdots & \vdots & \vdots & \vdots & \ddots
\end{bmatrix},
\end{equation}
and
\[
{\bf h} = [\,\cdots,\; -2ig_{2},\; -ig_{1},\; ig_{-1},\; 2ig_{-2},\;\cdots]^{\mathsf{T}}
\]
respectively.
Numerical simulations suggest that the finite-section approximation of the linearization in \eqref{CarlemanFourier.eq0} may exhibit exponential convergence; see Figure \ref{fig:kuramoto}. Nonetheless, a rigorous mathematical justification remains elusive, primarily because the matrix
\({\bf H}\) in \eqref{CarlemanFourier.eq0} lacks the upper triangular structure, which is  essential to the convergence analysis in \eqref{Carleman.exponential}. Furthermore, the intricate coupling within \({\bf H}\) poses significant challenges to establishing theoretical guarantees.

These challenges underscore the necessity of re-evaluating the linearization strategy for the dynamical system in \eqref{dynamicsystem.def}. In response, we propose a {\bf novel} framework that harnesses the periodicity of \(g(x)\) to construct a block-upper triangular state matrix while promoting sparsity within that structure. This revised approach retains the fundamental benefits of traditional Carleman techniques, such as facilitating rigorous convergence analyses, yet circumvents the structural shortcomings associated with direct Fourier-based linearization. Concretely, we introduce an extended state vector \(\tilde {\bf x} = [x_1, x_2]^{\mathsf{T}}\), where \(x_1 = x\) and \(x_2 = -x\). By \eqref{dynamicsystem.def}, \(\tilde {\bf x}\) satisfies
\begin{equation}\label{extendeddynamicsystem.def}
\dot {\tilde {\textbf{x}}} \;=\; \tilde{\bf g}(\tilde{\bf x}) \;:=\;
\begin{bmatrix} g_0 \\[6pt] -g_0 \end{bmatrix}
\;+\; \sum_{m=1}^\infty \begin{bmatrix} g_m & g_{-m} \\[4pt] -g_m & -g_{-m} \end{bmatrix}
e^{im\tilde{\bf x}},
\end{equation}
where \( e^{im\tilde{\bf x}} = [\,e^{imx_1},\; e^{imx_2}\,]^{\mathsf{T}},  m \ge 1 \). Unlike the original nonlinear system in \eqref{dynamicsystem.def}, this extended formulation ensures that \(\tilde{\bf g}(\tilde{\bf x})\) is a periodic function involving  {\bf nonnegative}  frequencies only, thus underpinning a more flexible and sparse Carleman-based representation. Through this construction, we reconcile the periodic nature of \(g(x)\) with the need for a tractable and sparse upper triangular (or block-upper triangular) matrix, ultimately yielding a more robust and scalable linearization framework.

\medskip

Set $z_1=e^{ix_1}$ and $z_2=e^{ix_2}$.   For  nonnegative integers $k_1$ and $k_2$ with $k_1+k_2\ge 1$,
we obtain from \eqref {extendeddynamicsystem.def} that
\begin{eqnarray}\label{derivative.wk-1}
\frac{d (z_1^{k_1} z_2^{k_2})}{dt} & \hskip-0.08in = & \hskip-0.08in  i(k_1-k_2) z_1^{k_1} z_2^{k_2} \Big(g_0+\sum_{m=1}^\infty g_m z_m^l+ g_{-m} z_2^m\Big)\nonumber\\
& \hskip-0.08in  = & \hskip-0.08in  i(k_1-k_2)\Big\{ g_0 z_1^{k_1} z_2^{k_2}+ \sum_{m=1}^\infty g_m z_1^{k_1+m} z_2^{k_2} + \sum_{m=1}^\infty  g_{-m} z_1^{k_1} z_2^{k_2+m}\Big\}.
\end{eqnarray}
Define the vector of all $k$-th order Fourier basis functions as
 $${\bf y}_k= \big[z_1^{k}, z_1^{k-1} z_2, \ldots, z_1 z_2^{k-1}, z_2^k\big]^{\mathsf{T}}\ \ {\rm  and} \ \
{\bf y}_k^0= \big[e^{ikx_0}, e^{i(k-2)x_0}, \ldots, e^{-i(k-2)x_0}, e^{-ikx_0}\big]^{\mathsf{T}} $$
for $k\ge 1$.
Then grouping all monomials $z_1^{k_1}z_2^{k_2}$ with $k_1+k_2=k$ together and applying \eqref{derivative.wk-1},
we have
\begin{equation}\label{derivative.wk}
\dot {\bf y}_k= \sum_{l=k}^\infty {\bf B}_{k,l} {\bf y}_l\ \ {\rm for} \ k\ge 1,
\end{equation}
 where  for every $1\le k\le l$, matrix ${\bf B}_{k, l}=\big[b_{k, l; p, q}\big]_{0\le p\le k, 0\le q\le l}$ 
depends on Fourier coefficients $g_{\pm (l-k)}$ in \eqref{Fourierexpansion.def} and it is  given by
\begin{equation*}\label{Bkl.def}
{ b}_{k,l; p, q}=\left\{\begin{array}
    {ll}
    i(k-2p) g_{l-k} & {\rm if} \ q=p\\
    i (k-2p) g_{k-l} & {\rm if} \ q=l-k+p\\
0 & {\rm otherwise}.
    \end{array}\right.
\end{equation*}
By combining all ODEs  in \eqref{derivative.wk} and defining a new state vector $${\bf y}=\big[{\bf y}_1^{\mathsf{T}}, {\bf y}_2^{\mathsf{T}},\ldots, {\bf y}_N^{\mathsf{T}}, \ldots \big]^{\mathsf{T}},$$ we obtain the following infinite-dimensional linear system
\begin{equation}
\label{CarlemanFourier.def}
\dot {\bf y}={\bf B} {\bf y} =
\begin{bmatrix}
 {\bf  B}_{1,1} & {\bf   B}_{1,2} & \dots &  {\bf   B}_{1, N}  & \cdots\\
      & {\bf B}_{2,2} & \cdots & {\bf   B}_{2, N} & \cdots\\
      & & \ddots & \vdots & \ddots\\
      & & & {\bf   B}_{N,N}& \ddots \\
      & & & & \ddots
\end{bmatrix} \begin{bmatrix} {\bf y}_1 \\ {\bf y}_2 \\ \vspace{0.07cm} \vdots \vspace{0.07cm} \\ {\bf y}_N \\ \vdots
\end{bmatrix}
\end{equation}
with initial condition
$$ {\bf y}(0)=\big[({\bf y}_1^0)^{\mathsf{T}},\ ({\bf y}_2^0)^{\mathsf{T}},\ \ldots\big]^{\mathsf{T}}.$$
Recall that
 the state matrix ${\bf A}$
 in the classical Carleman linearization
 exhibits an upper-triangular structure. Similarly, we observe that the state matrix ${\bf B}$ in \eqref{CarlemanFourier.def} forms a block-upper triangular matrix, and  its diagonal blocks are {\bf diagonal} matrices depending only on $g_0$. Based on this observation, we refer to the resulting linearization in \eqref{CarlemanFourier.def} as the {\em Carleman-Fourier linearization} of the nonlinear dynamical system detailed in \eqref{dynamicsystem.def}.

Analogously to the state matrix ${\bf H}$ in \eqref{CarlemanFourier.eq0}, which is {\bf not} a bounded operator on the space of square-summable sequences, the state matrix ${\bf B}$ in \eqref{CarlemanFourier.def} shares a similar characteristic. To address this challenge, we introduce a {\em finite-section approximation} of the Carleman-Fourier linearization in \eqref{CarlemanFourier.def}, which provides a 
computationally efficient approach, by
\begin{equation}
\label{Carleman.eq8}
\begin{bmatrix}
    {\Dot {\bf   y}}_{1,N} \\ {\Dot{ {\bf   y}}}_{2,N} \\ \vdots \\ {\Dot{{\bf   y}}}_{N,N}
\end{bmatrix}
    =\begin{bmatrix}
     {\bf   B}_{1,1}& {\bf B}_{1, 2} & \ldots & {\bf   B}_{1, N} \\
      &{\bf  B}_{2,2}& \ldots & {\bf  B}_{2, N} \\
      & & \ddots & \vdots \\
      & & & {\bf  B}_{N,N}
      \end{bmatrix}
      \begin{bmatrix}
    {\bf  y}_{1,N} \\  {\bf  y}_{2,N}\\ \vdots \\ {\bf  y}_{N,N}
\end{bmatrix}
\end{equation}
with the initial  condition
${\bf  y}_{k,N}(0)={\bf  y}_k^0$ for   $1\le k\le N$.
  In this paper, we show that
${\bf y}_{1, N}$ converges exponentially to ${\bf y}_1$  in a quantifiable time range.

\begin{theorem}\label{main.thm}
Consider a periodic function $g$ as defined in \eqref{Fourierexpansion.def} and fulfilling the condition detailed in \eqref{fouriercoefficient.decay}. Let $x(t)$ represent the state of the dynamical system given by \eqref{dynamicsystem.def}, and ${\bf y}_{1, N}(t)$ as specified in \eqref{Carleman.eq8}. Let us  write   ${\bf y}_{1, N}(t)=[y_{1, N}^+(t), y_{1, N}^-(t)]^{\mathsf{T}}$ and define
$$T_0= \frac{1}{2D}\Big(\frac{1}{ \sqrt{r}}-1\Big)^2,$$
where constants $D>0$ and $r\in (0, 1)$ are from \eqref{fouriercoefficient.decay}. Then, for $N\ge 1$, the inequality
\begin{equation}\label{mainestimate}
\big|y_{1, N}^\pm (t)-e^{\pm ix(t)}\big|\le \frac{\sqrt{2D t}}{N(1-r)} \hspace{0.05cm} \big(1+\sqrt{2D t}\big)^{2N} r^N
\end{equation}
hold true for all $0 \leq t \le  T_0$.
\end{theorem}

The result in Theorem \ref{main.thm} elucidates the direct relationship between the solution of the approximate Carleman-Fourier linearization, as denoted in \eqref{Carleman.eq8}, and the complex exponential of the solution of the original nonlinear system given by \eqref{dynamicsystem.def}. The detailed proof of Theorem \ref{main.thm} is presented in Section \ref{proof.section}.

We select  an integer \( N_0 \ge 1 \) such that
\begin{equation}\label{N0.def}
N_0\ge \frac{2\sqrt{2D T_0}}{1-r}.
\end{equation}
Consequently, for every \( N \ge N_0 \), the functions \( y_{1, N}^\pm (t) \) are non-zero and continuous over the interval
\([0, T_0]\), as established in \eqref{mainestimate}. These functions can be expressed as:
\begin{equation} \label{thetaN.def0}
y_{1, N}^\pm (t)= |y_{1, N}^\pm(t)| \ e^{i\vartheta_{1, N}^\pm (t)}\ \ {\rm with} \ \ \vartheta_{1, N}^\pm (0) = \pm x_{0}.
\end{equation}
It is noteworthy that for any complex number \( z \) satisfying \( |e^{iz}-1|\le \epsilon\le 1/2 \), the condition \( |z \ {\rm mod} \ 2\pi|\le  4 \epsilon \) holds. This insight, combined with the conclusions from Theorem \ref{main.thm}, indicates that the functions \( \vartheta^\pm_{1, N}(t)  \), \( N \ge 1 \), as defined in \eqref{thetaN.def0}, offer an exponential estimate to the states \( \pm x(t) \) over the time interval \( 0\le t\le T_0 \).

\begin{corollary}\label{main.cor1}
In line with the assumptions outlined in Theorem \ref{main.thm}, we  select an integer $N_0$ that adheres to \eqref{N0.def}. We then define $\vartheta^\pm_{1, N}(t)$ for any $N\ge N_0$ as per \eqref{thetaN.def0}. Then, for all $N\ge N_0$ and $0 \leq t \leq T_0$, it follows that
\begin{equation*}\label{mainestimate.cor.eq}
|\vartheta_{1, N}^\pm (t)\mp x(t)|\le \frac{4\sqrt{2D t}}{N(1-r)} \big( \big(1+\sqrt{2D t}\big)^2 r\big)^N.
\end{equation*}
\end{corollary}

\begin{remark}\label{remark3}
{\rm From 
Theorem \ref{main.thm} and Corollary \ref{main.cor1},
we note that our estimates to the convergence rate and the duration of exponential convergence of the finite-section approximation to the Carleman-Fourier linearization are {\bf independent} of the initial condition $x_0$ and the value $g(0)$. This global convergence property for the Carleman-Fourier linearization contrasts with the traditional Carleman linearization approach, where the convergence of its finite-section approximations is highly dependent on $x_0$ and $g(0)$, particularly in terms of proximity to the equilibrium point zero, see Figure \ref{fig:kuramoto} for numerical demonstration.
}\end{remark}

\section{Carleman-Fourier linearization of systems with multiple fundamental frequencies}
\label{CarlemanFourierGeneral.section}

Dynamical systems with quasi-periodic vector fields possess a remarkable richness and complexity
in their internal dynamics
\cite{ Akiba2023, Amini2022,  Brunton2016, Liu2021}.
 In this section, we tailor our Carleman-Fourier linearization method proposed in the last section to the following family of $d$-dimensional dynamical systems with quasi-periodic vector field,
\begin{equation}\label{dynamicsystem.multiple}
\dot{\bf x}(t) = {\bf g}({\bf x}(t))
\end{equation}
 with the initial condition ${\bf x}(0)= {\bf x}_0$, where ${\bf x}=[x_1, \ldots, x_d]^{\mathsf{T}}\in {\mathbb R}^d$ and ${\bf g}({\bf x})=[g_1({\bf x}), \ldots, g_d({\bf x})]^{\mathsf{T}}$. The functions
 \begin{equation} \label{dynamicsystem.eqb}
{g}_p( {\bf x})=\hskip-0.08in
\sum_{{\pmb \alpha}_1, \ldots, {\pmb \alpha}_L\in {\mathbb Z}^d} {g}_{p; {\pmb \alpha}_1, \ldots, {\pmb \alpha}_L}  \hspace{0.05cm} e^{ i
(\tau_1 {\pmb \alpha}_1 +\cdots+\tau_L {\pmb \alpha}_L)^{\mathsf{T}}{\bf x}}, \ 1\le p\le d,
\end{equation}
 consist of real-valued quasi-periodic functions with multiple distinct fundamental frequencies $\tau_1, \ldots, \tau_L >0$, where $\mathbb{Z}^d$ is the set of $d$-dimensional vectors with integer components.
 These functions are characterized by Fourier coefficients that exhibit exponential decay. Specifically, there exist positive constants $D>0$ and $r\in (0, 1)$ ensuring that the following decay condition is satisfied for every $k\ge 0$,
\begin{equation} \label {dynamicsystem.eqbb}
\sup_{1\le p\le d} \sum_{|{\pmb \alpha}_1|+\ldots+|{\pmb \alpha}_L|=k}
|g_{p; {\pmb \alpha}_1, \ldots, {\pmb \alpha}_L}|\le \frac{2^dD}{\tau_1+\cdots+\tau_L} \hspace{0.05cm} r^k.
    \end{equation}


It is evident that our illustrative dynamical system described in \eqref{dynamicsystem.def} represents a specific case of the aforementioned dynamical system in a one-dimensional setting with a single fundamental frequency.  Our illustrative Kuramoto model of coupled oscillators, as formulated in \eqref{Kuramoto.def}, exemplifies a dynamical system with a single fundamental frequency $\tau_1=1$. In this model, the constants in \eqref{dynamicsystem.eqbb} are defined as $D= 2^{-d-1} K r^{-2}$ and $r\in (0, 1)$.

We define the extended state vector as
\[
\tilde {\bf x}:=[\tau_1 {\bf x}^{\mathsf{T}}, \ldots, \tau_L {\bf x}^{\mathsf{T}},
-\tau_1 {\bf x}^{\mathsf{T}}, \ldots, -\tau_L {\bf x}^{\mathsf{T}}]^{\mathsf{T}}\in {\mathbb R}^{2dL}
\]
and the initial extended state vector as
\[
\tilde {\bf x}_0:= [\tau_1 {\bf   x}_0^{\mathsf{T}},  \ldots, \tau_L {\bf   x}_0^{\mathsf{T}},
-\tau_1 {\bf   x}_0^{\mathsf{T}}, \ldots, -\tau_L {\bf   x}_0^{\mathsf{T}}]^{\mathsf{T}}\in {\mathbb R}^{2dL}.
\]
 Subsequent verification confirms that the extended state vector $\tilde {\bf x}$ satisfies the following nonlinear dynamical system
\begin{equation*}  
{\Dot{\tilde {\bf   x}}}(t)= {\bf f}(\tilde {\bf x}(t))
\ \ {\rm for } \  t\ge 0,
\end{equation*}
 with the initial condition set as $\tilde {\bf x}(0)=\tilde {\bf x}_0$, where
 $\mathbb{Z}_+^{2dL}$ is the set of $2dL$-dimensional vectors with nonnegative integer components, and the function
\[
{\bf f}(\tilde {\bf x})=\sum_{{\pmb \gamma}\in {\mathbb Z}_+^{2dL}}
[   f_{1; {\pmb \gamma}}, \Sp \ldots \Sp,  f_{2dL; {\pmb \gamma}}]^{\mathsf{T}}  e^{i {\pmb \gamma}^{\mathsf{T}} \tilde  {\bf   x}}
\]
is periodic with respect to the extended variable $\tilde {\bf x}$ and features nonnegative frequencies, in which
$
{\pmb \gamma}=\big[({\pmb \alpha}_{1})_+^{\mathsf{T}}, \Sp \ldots \Sp, ({\pmb \alpha}_{L})_+^{\mathsf{T}},
({\pmb \alpha}_{1})_-^{\mathsf{T}}, \Sp  \ldots \Sp, ({\pmb \alpha}_{L})_-^{\mathsf{T}} \big]^{\mathsf{T}}
$
for certain ${\pmb \alpha}_1, \ldots, {\pmb \alpha}_L\in {\mathbb Z}^d$. Here for ${\pmb \alpha}\in {\mathbb Z}^d$ we let \(\pmb{\alpha}_+\) be the vector \([\max(\alpha_1, 0), \ldots, \max(\alpha_d, 0)]^{\mathsf{T}}\in {\mathbb Z}_+^d\) and \(\pmb{\alpha}_-\) be the difference \(\pmb{\alpha}_+ - \pmb{\alpha}=[-\min(\alpha_1, 0), \ldots, -\min(\alpha_d, 0)]^{\mathsf{T}}\in {\mathbb Z}_+^d\).
The function $f_{j; {\pmb \gamma}}(t)$ equals zero except when it takes the form $f_{j; {\pmb \gamma}}(t)= (-1)^m \Sp \tau_l \Sp g_{p; {\pmb \alpha}_1, \ldots, {\pmb \alpha}_L}(t)$, applicable for $j= m Ld+ (l-1) d+p$ for specific integer values of $1\le l\le L$, $1\le p\le d$, and $0\le m\le 1$. 

 Given a vector \(\pmb{\alpha} = [\alpha_1, \ldots, \alpha_d]^{\mathsf{T}} \in \mathbb{Z}^d\) and an integer $k\ge 1$, we define the norm \(|\pmb{\alpha}|\) as the sum of the absolute values of its components: \(|\alpha_1| + \cdots + |\alpha_d|\), and
 denote   the subset of \(\mathbb{Z}_+^d\) consisting of vectors \(\pmb{\alpha}\) with the norm \(|\pmb{\alpha}|\) equaling to \(k\) by $\mathbb{Z}_{+, k}^d$.
For every $k\ge 1$, we denote all the
$k$-th order terms as
\[
{\bf z}_k= \big[e^{i{\pmb \gamma}^{\mathsf{T}}\tilde {\bf x}}\big]_{{\pmb \gamma}\in {\mathbb Z}^{2dL}_{+, k}}
\]
and their corresponding initial conditions by
\[
{\bf z}_{k}^0= \big[e^{i{\pmb \gamma}^{\mathsf{T}}\tilde {\bf x}_0}\big]_{{\pmb \gamma}\in {\mathbb Z}^{2dL}_{+, k}},
\]
and  for $1\le k\le l$, we define the block matrices
\begin{equation}\label{function.multiple.eq5+}
    {\bf   F}_{k,l}=\Big[\Sp i\sum_{j=1}^{2dL} \Sp \gamma_j \Sp f_{j; {\pmb \delta}-{\pmb \gamma}}\Sp \Big]_{{\pmb \gamma}\in {\mathbb Z}_{+,k}^{2dL}, \Sp {\pmb \delta}\in {\mathbb Z}_{+,l}^{2dL}}
\end{equation}
where ${\pmb \gamma}=[\gamma_1, \ldots, \gamma_{2dL}]^{\mathsf{T}}$. Following the methodology outlined in \eqref{CarlemanFourier.def}, we define the Carleman-Fourier linearization of the dynamical system in \eqref {dynamicsystem.multiple} by
\begin{equation} \label{Carleman.multiple.eq7}
\begin{bmatrix}
   \dot{ { {\bf  z}}}_{1}\\
    \dot{  {\bf z}}_{2} \\
      \vdots \\
       \dot{{\bf  z}}_{N} \\
       \vdots
\end{bmatrix}
=
\begin{bmatrix}
   {\bf  F}_{1,1} &  {\bf   F}_{1,2}& \dots  &  {\bf   F}_{1, N}  & \cdots\\
      & {\bf F}_{2,2} & \cdots & {\bf   F}_{2, N} & \cdots\\
      & & \ddots & \vdots & \ddots\\
      & & & {\bf   F}_{N,N} & \cdots \\
      & & & & \ddots
\end{bmatrix}
\begin{bmatrix}
     {\bf  z}_{1}\\
     {\bf z}_{2} \\
      \vdots \\
       {\bf  z}_{N}  \\
       \vdots
\end{bmatrix}
\end{equation}
 with  the initial $ {\bf z}_k(0)=  {\bf z}_k^0$ for all $k\ge 1$,
and its finite-section approximation of order $N\ge 1$ by
\begin{equation}
\label{Carleman.multiple.eq8}
\begin{bmatrix}
    {\Dot {{\bf   z}}}_{1,N} \\ {\Dot{ {\bf   z}}}_{2,N} \\ \vdots \\ {\Dot{{\bf   z}}}_{N,N}
\end{bmatrix}
    =\begin{bmatrix}
     {\bf   F}_{1,1} &   {\bf   F}_{1,2} & \ldots &  {\bf   F}_{1, N}  \\
      &{\bf  F}_{2,2} & \ldots & {\bf  F}_{2, N} \\
      & & \ddots & \vdots \\
      & & & {\bf  F}_{N,N}
      \end{bmatrix}
      \begin{bmatrix}
    {\bf  z}_{1,N} \\  {\bf  z}_{2,N} \\ \vdots \\ {\bf  z}_{N,N}
\end{bmatrix}
\end{equation}
 with the initial $ {\bf z}_{k,N}(0)=  {\bf z}_k^0$ for $1\le k\le N$.

 Set $ {\pmb \xi}_{k,N}= {\bf z}_{k, N}- {\bf z}_{k}$. As  ${\bf F}_{k, k}, k\ge 1$, are diagonal matrices with diagonal
entries being pure imagery, we can show that
\begin{eqnarray*}
\|{\pmb \xi}_{k, N}(t)\|_\infty
 &  \le &   2^d D k \int_0^t \sum_{l=k+1}^N r^{l-k} \ \| {\pmb \xi}_{l, N}(s)\|_\infty \ ds
  + \frac{2^d D kr}{1-r} r^{N-K} t \ \ {\rm for} \ t\ge 0
\end{eqnarray*}
hold for all $1\le k\le N$.
Then, following the argument used in the proof of Theorem \ref{main.thm}, we can show  exponential convergence of the primary block $ {\bf z}_{1, N}, N\ge 1$,
of the finite-section approach in  \eqref{Carleman.multiple.eq8}.

\begin{theorem}\label{multiplemain.thm} Consider ${\bf g}$ as a quasi-periodic function with fundamental frequencies $\tau_1, \ldots, \tau_L$, satisfying \eqref{dynamicsystem.eqb} and \eqref{dynamicsystem.eqbb}. Let ${\bf x}(t)=[x_1(t), \ldots, x_d(t)]^{\mathsf{T}}$ be the state vector of the nonlinear system in \eqref{dynamicsystem.multiple}, and  ${\bf z}_{1, N}(t)=[ {z_{1; 1, N}}(t),\ldots,  {z}_{2dL; 1, N}(t)]^{\mathsf{T}}$ represent the first block in \eqref{Carleman.multiple.eq8}. Assuming $$T_1=\frac{(1-r)^2}{2^dD r^2}$$ and selecting an integer $N_0 >0$ such that
$$2(1-r)^{-1} \sqrt{2^dD T_1}   \le  N_0,$$
with $D>0$ and $r\in (0, 1)$ as constants in \eqref{dynamicsystem.eqbb}, we express each state variable
${{z}_{j; 1, N}}(t)= |{z_{j; 1, N}}(t)| \hspace{0.05cm} e^{i \vartheta_{j; 1, N}(t)}$
for some real-valued continuous function $\vartheta_{j; 1, N}(t)$ with initial  $[\vartheta_{1; 1, N}(0), \ldots, \vartheta_{2dL; 1, N}(0)]^{\mathsf{T}} = \tilde {\bf x}_0$. Then, for every $N\ge N_0$ and $0\le t\le T_1$, it holds that
\begin{eqnarray*} 
  &   &   \big| {\vartheta}_{p+d(l-1)+mdL,N}(t) -(-1)^m\tau_lx_{p}(t)\big|  
 \le  \frac{4 \sqrt{2^dD t}}{N(1-r)} \big( \big(1+\sqrt{2^dD t}\big)^2 r\big)^N,
\end{eqnarray*}
where $m\in \{0, 1\}, 1\le l\le L$ and $1\le p\le d$.
\end{theorem}

\section{Applications: Kuramoto model}
The renowned first-order Kuramoto model is  expressed by the following equation:
\begin{equation}\label{Kuramoto.def}
\dot{\theta}_{p} = \omega_{p} + \frac{K}{d} \sum_{q=1}^{d} \sin(\theta_{q} - \theta_{p})~~\textrm{for}~~  1 \leq p \leq d,
\end{equation}
where $\theta_p$ denotes the phase of the
$p$-th oscillator with natural frequency $\omega_p$ for $1 \leq p \leq d$, and $K \neq 0$ represents the coupling strength between oscillators. This model serves as a quintessential example of dynamical systems
and is foundational for studying nonlinear dynamical systems, particularly in the context of coupled oscillators. This model offers valuable insights into synchronization phenomena observed in natural and technological systems, capturing how individual components, despite differing intrinsic frequencies, can achieve collective coherence through mutual interactions \cite{Acebron2005, bronski2021, dietert2016, guo2021, heggli2019, ji2014, kuramoto1984}.

\begin{figure*}[ht] 
\centering
         \includegraphics[width=5.8cm, height=4.8cm]{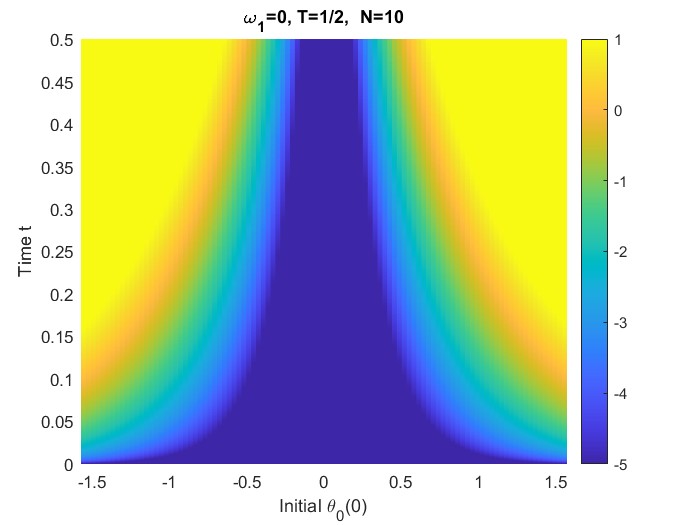}
     \includegraphics[width=5.8cm, height=4.8cm]
    {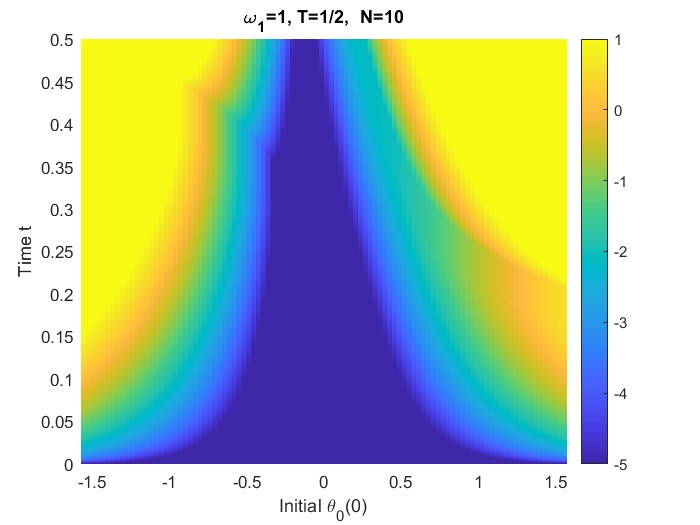}
    \\
                 \includegraphics[width=5.8cm,  height=4.8cm]
     {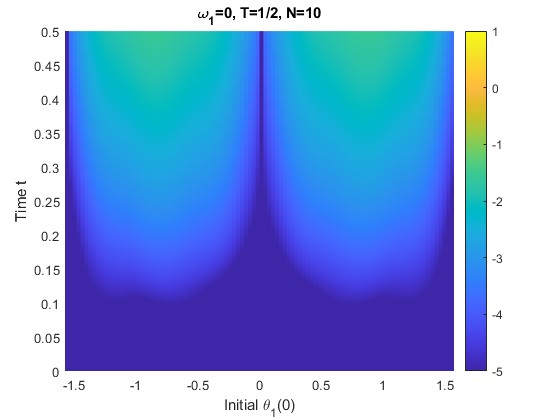}
   \includegraphics[width=5.8cm,  height=4.8cm]      {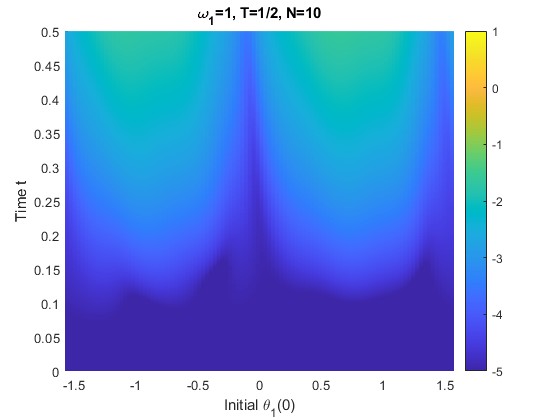}
       \\
            \includegraphics[width=5.8cm,  height=4.8cm]
     {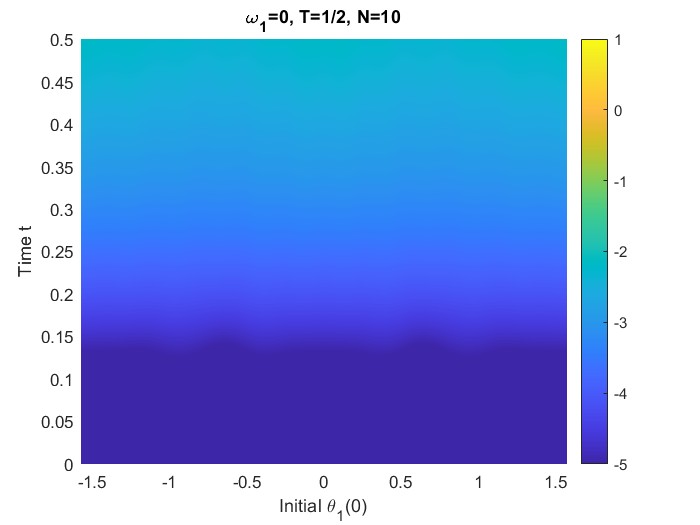}
     \includegraphics[width=5.8cm,  height=4.8cm]
     {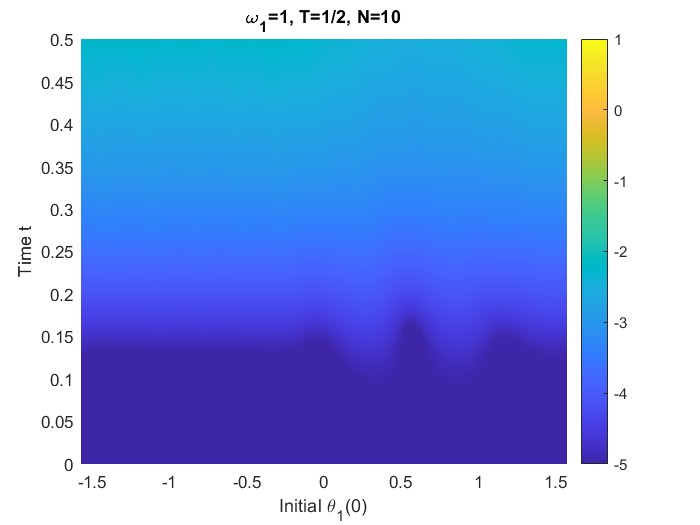}

\caption{Plotted  from top to bottom are
the approximation error
$E_{\rm C}(\theta_1(0), \omega_1, N, t)$
in \eqref{ECNt.def} for the  finite-section approach to the classical Carleman linearization,
the approximation error $E_{\rm CF}(\theta_1(0), \omega_1,  N, t)$
in \eqref{ECFNt.def} for the finite-section approach to the Carleman-Fourier linearization,
and the approximation error $ E(\theta_1(0), \omega_1, N, t)$ in \eqref{ENt.def}
 for the  finite-section approach to the conventional linearization   respectively,
  where $-\pi/2\le \theta_1(0)\le \pi/2, 0\le t\le 1/2$.
Here $\tilde K=1, N=10$ and $\omega_1=0$ (all left) and  $\omega_1=1$ (all right).
  }
  \label{fig:kuramoto}
\end{figure*}

In this paper, we present numerical simulations to demonstrate the effectiveness of the finite-section approach applied to both the classical Carleman linearization and the proposed Carleman-Fourier linearization for the Kuramoto model of coupled oscillators. 
As shown in Figure \ref{fig:kuramoto}, we observe that the finite-section approach to the proposed Carleman-Fourier linearization exhibits significantly better approximation properties compared to the classical Carleman linearization, particularly when the initial phase and natural frequency are away from zero. Unlike the classical approach, the proposed method imposes no restrictions on the natural frequency and the initial phase, as noted in Remark \ref{remark3}.

In this paper, we normalize the Kuramoto model so that
\vspace{-0.2cm}
\begin{equation}
\sum_{p=1}^d \theta_p(0) = \sum_{p=1}^d \omega_p = 0\ \ {\rm and}\ \ |K|=d,
\end{equation}
otherwise, replace the phases $\theta_p$ and
natural frequencies $\omega_p$ by  alternative phases $\tilde \theta_p$ and frequencies
$\tilde \omega_p, 1\le p\le d$,
which are given by
$$\tilde \theta_p(t)= \theta_p \left(\frac{d}{|K|}t\right)- \frac{ \sum_{q=1}^d w_q}{|K|}
t- \frac{1}{d}\sum_{q=1}^d \theta_q(0)$$
and
$$ \tilde \omega_p=  |K|^{-1} \Big(d\omega_p-  \sum_{q=1}^d w_q\Big)\ \
 {\rm for} \ 1\le p\le d.$$

For the case where $d = 2$ and with the normalization of initial phases and natural frequencies, the phase $\theta_1$ of the oscillator adheres to the dynamical system in the form described in \eqref{dynamicsystem.def},
\begin{equation}\label{kuramoto.def.one+}
\dot{\theta}_1(t) = \omega_1 + \tilde K \sin(2\theta_1(t))\ \ {\rm for}\  t \ge 0,
\end{equation}
where $\tilde K=-K/d=\pm 1$.
One may verify that
the state variable
$\theta_1$ diverges when $|\omega_1|>1$,
and it
converges to one of the equilibria  $\{-
\frac{\tilde K}{2} \arcsin \omega_1, \frac{\tilde K}{2} \arcsin \omega_1- \frac{\pi}{2}\}+ \pi {\mathbb Z}$ when $|\omega_1|\le 1$.

   For the nonlinear dynamical system in \eqref{kuramoto.def.one+},
   the linear system of size $N\ge 1$ associated with the finite-section approximation to its Carleman linearization is given by
   \begin{equation}
\label{CarlemanKuramoto.def}
\begin{bmatrix} \dot \theta_{1, N} \\  \dot \theta_{2, N} \\ \vspace{0.07cm} \vdots \vspace{0.07cm} \\ \dot \theta_{N, N}
\end{bmatrix} =
\begin{bmatrix}
a_{11} & a_{12} &    \dots &  a_{1N}   \\
  \omega_1    & a_{22} &  \cdots &
a_{2N} \\
      &\ddots & \ddots & \vdots \\
      & &\omega_1 &   a_{NN}
\end{bmatrix} \begin{bmatrix} \theta_{1, N} \\ \theta_{2, N} \\ \vspace{0.07cm} \vdots \vspace{0.07cm} \\ \theta_{N,N}
\end{bmatrix}+
\begin{bmatrix} \omega_1 \\  0 \\ \vspace{0.07cm} \vdots \vspace{0.07cm} \\ 0
\end{bmatrix},
\end{equation}
where
$$a_{nn'}= \left\{\begin{array}{ll} \hskip-.08in
     \frac{\tilde K n 2^{n'-n+1} (-1)^{(n'-n)/2}}{(n'-n+1)!} & \hskip-.08in {\rm if}\ \
    n'-n\in 2{\mathbb Z}\cap [0, N-1]\\
 \hskip-0.08in     0 &  \hskip-.08in {\rm otherwise}.
     \end{array}
     \right.$$
As shown in \eqref{Carleman.exponential}, the first variable $\theta_{1, N}$,  $N\ge 1$, in \eqref{CarlemanKuramoto.def}  suggests exponential convergence to the original phase $\theta_1$ in a short time range
when the initial state $\theta_1(0)$ is not far from zero and the natural frequency $\omega_1$ takes zero value. Numerical confirmation for this behavior is presented in the top plot of Fig. \ref{fig:kuramoto}, where
\begin{equation}\label{ECNt.def}
  E_{\rm C}(\theta_1(0), \omega_1,  N, t)   =   \sup_{0\le s\le t} \log_{10} \min \Big\{ 10,
\max \big\{|\theta_{1, N}(s)-\theta_1(s)|, 10^{-5}\big\}\Big\}
\end{equation}
for $0\le t\le T$, denotes the approximation error  for the finite-section  approach to the classical Carleman linearization of the Kuramoto model in \eqref{kuramoto.def.one+} on the time interval $[0, T]$ in the logarithmic scale.  This reconfirms that the classical Carleman linearization is a prominent method
 to linearize a dynamical system in the neighborhood of the origin.
We believe that the finite-section method to the  classical Carleman linearization exhibits exponential convergence even when the equilibrium point requirement $g(0) = 0$ in \eqref{Carleman.requirement} is not met. Numerical simulations (see top right plot of Fig. \ref{fig:kuramoto}) support this conjecture for the Kuramoto model in \eqref{Kuramoto.def} when the initial phase $\theta_1(0)$ is near zero, where $g(0) = \omega_1=1$. However, comparing with the case that $\omega_1=0$, the finite-section approximation  converges to the original phase in a shorter time range and  with slow convergence.

\medskip

For the Kuramoto model in \eqref{kuramoto.def.one+}, we define the extended phase $\tilde{\pmb \theta} = [\theta_1, -\theta_1]^{\mathsf{T}}$.
We observe that the entries of matrices ${\bf B}_{k, l} = \big[b_{k, l; p, q}\big]_{0\le p\le k, 0\le q\le l}$ in the Carleman-Fourier linearization are zero except for ${b}_{k, k; p, p} = i(k-2p) \omega_1$, ${b}_{k, k+2; p, p} = (k-2p) \tilde  K/2$, and ${b}_{k, k+2; p, p+2} = -(k-2p) \tilde K/2$ for $0 \leq p \leq k$ and $k \geq 1$. Let
$${\bf v}_{1, N}(t) = [v_{1,N}^+(t), v_{1,N}^-(t)]^{\mathsf{T}}$$
be the first block in the finite-section approximation of the Carleman-Fourier linearization associated with the Kuramoto model of coupled oscillators.
Take arbitrary $T^*\in (0, 1)$. By Theorem \ref{main.thm}, we conclude that $v_{1,  N}^\pm (t), N\ge 1$, converge exponentially to $e^{\pm i\theta_1(t)}$ in the time range $[0, T^*]$, {\bf no matter what} the initial phase $\theta_1(0)$ and the natural frequency $\omega_1$ are located, and the exponential convergence rate  could be bounded by
a positive number $r\in (0,1)$ independent of the initial $\theta_1(0)$ and the natural frequency $\omega_1$. Shown on the middle plots  of Figure \ref{fig:kuramoto} are the approximation error
\begin{eqnarray}\label{ECFNt.def}
  &    & \hskip-0.4in  E_{\rm CF}(\theta_1(0), \omega_1, N, t)     =   \sup_{0\le s\le t}   \log_{10} \min\Big\{ 10,
\nonumber \\
 &  & \hskip-0.40in \qquad\quad
  \max \big\{|v_{1, N}^+(s)- e^{i\theta_1(s)}|,  |v_{1, N}^-(s)-e^{-i\theta_1(s)}|, 10^{-5}\big\}  \Big\}
 \end{eqnarray}
of the finite-section approach to the Carleman-Fourier linearization of the Kuramoto model in \eqref{kuramoto.def.one+} in the logarithmic scale, where
$0\le t\le 1/2$ and $-\pi/2\le \theta_1(0)\le \pi/2$.
As
$$\max \big\{|v_{1, N}^+(t)- e^{i\theta_1(t)}|,  |v_{1, N}^-(t)-e^{-i\theta_1(t)}|\big\}=
\max \big\{|v_{1, N}^+(t) e^{-i\theta_1(t)}-1|,  |v_{1, N}^-(t)e^{i\theta_1(t)}-1|\big\},$$  we conclude that
 the imaginary part of the logarithm of the first block ${\bf v}_{1, N}(t)$ of the state vector in the finite-section approach converges exponentially to the extended state $\tilde{\pmb \theta}(t)$ in the time range $[0, 1)$, cf. Corollary \ref{main.cor1}.
  It is observed from the plots on the middle left and right that the
  approximation errors $ E_{\rm CF}(\theta_1(0), \omega_1, N, t)$ are quite similar for $\omega_1=0$ and $\omega_1=1$, except when $\theta_1(0)$ is close to zero and $\pm \pi/2$.
Compared with classical Carleman linearization,
the proposed Carleman-Fourier linearization technique exhibits much better approximation except that the natural frequency $\omega_1$ and the initial $\theta_1(0)$ are very close to the origin.
This well-approximation phenomenon can be observed from  Fig. \ref{fig:kuramoto} for $\omega_1=0$ and $1$.
We  also test the performance of the classical Carleman linearization and the Carleman-Fourier linearization for the  Kuramoto model with $d=3$ in our subsequent paper \cite{Chen2024}. It is also noticed that
its Carleman-Fourier linearization
delivers more accurate linearizations over more extensive range of the initial phases, and outperforms the classical Carleman
linearization when the initial phases of oscillators are not close to zero.

\smallskip

  For the Kuramoto model in \eqref{kuramoto.def.one+}, the Fourier coefficients $g_n, n\in {\mathbb Z}$,
of the governing function  takes zero values except $g_0=\omega_1$ and $g_{\pm 2}= \pm  \tilde K/(2i) $. Then
 its conventional linearization  in \eqref{CarlemanFourier.eq0}
has  the state matrix ${\bf H}=[in g_{n'-n}]_{n, n'\in {\mathbb Z}_0}$ and nonhomogeneous term ${\bf h}=[in g_{-n}]_{n\in {\mathbb Z}_0}$
given by
\begin{equation*} 
{\bf H} =
\begin{bmatrix}
\ddots & \vdots&  \vdots & \vdots & \vdots & \vdots &  \vdots& \ddots\\
\cdots & -3i\omega_1 & 0 & -3\tilde K/2 & 0 & 0 &  0 & \cdots\\
\cdots & 0 & -2i\omega_1 & 0 & 0 & 0 & 0 & \cdots\\
\cdots & \tilde K/2 & 0 & -i\omega_1 & -\tilde K/2 &  0 &  0 & \cdots\\
\cdots & 0 & 0 & -\tilde K/2 & i\omega_1 & 0 & \tilde K/2&  \cdots\\
\cdots & 0 & 0 & 0 & 0 & 2i\omega_1 & 0 & \cdots\\
\cdots & 0  & 0 & 0 & -3\tilde K/2 & 0 & 3i\omega_1 & \cdots\\
\ddots & \vdots & \vdots & \vdots & \vdots & \ddots
\end{bmatrix},
\end{equation*}
and
$ {\bf h} = [\,\cdots, 0,  -\tilde K, 0,\; 0,\; -\tilde K,\; 0, \cdots]^{\mathsf{T}} $
respectively.
Following the finite-section procedure \eqref{Carleman.eq8},   for the linearization
\eqref{CarlemanFourier.eq0} of
Kuramoto model  \eqref{kuramoto.def.one+} we may define its  finite-section approximation of order $N\ge 2$
as follows:
\begin{equation}\label{conventional.remark.eq1}
\dot {\bf w}_N={\bf H}_N{\bf w}_N+{\bf h}_N
\end{equation}
with  initial ${\bf w}_N(0)=[e^{-iN\theta_1(0)}, \ldots,  e^{-i\theta_1(0)}, e^{i\theta_1(0)},  \ldots, e^{Ni\theta_1(0)}]^{\mathsf{T}}$,
where ${\bf w}_N=[w_{-N, N}$, $\ldots,    w_{-1, N},   w_{1, N},   \ldots, w_{N, N}]^{\mathsf{T}}$,
 ${\bf H}_N= [in g_{n'-n}]_{n, n'\in\{-N, \ldots, -1, 1, \ldots, N\}}$
 is a principal submatrix of the state matrix ${\bf H}$,
and ${\bf h}_N=
 [0,\cdots, 0,  -\tilde K, 0, 0, -\tilde K, 0,  \cdots, 0]^{\mathsf{T}}$.
Since the governing function  of the dynamical system
 \eqref{kuramoto.def.one+}
   does not contain odd frequency components,  the odd-indiced components of the state vector
 ${\bf w}_N$ satisfies the following  homogeneous linear system:
 \begin{eqnarray*}
 \begin{bmatrix}
 \dot w_{-\tilde N, N}\\
  \dot w_{-\tilde N+2, N}\\
  \vdots\\
  \vdots\\
 \dot w_{\tilde N-2, N}\\
 \dot w_{\tilde N, N}
\end{bmatrix} & = &
\begin{bmatrix}
-\tilde N &   \\
 & -\tilde N+2\\
 & & \ddots\\
 & & &  \tilde N-2\\
 & & & &  \tilde N
\end{bmatrix}\nonumber \\
& & \times
\begin{bmatrix}
i\omega_1 & \tilde K/2  \\
-\tilde K/2 & i\omega_1 &  \tilde K/2 \\
 & \ddots & \ddots & \ddots \\
& & \ddots & \ddots & \ddots \\
&& & \ddots & \ddots & \ddots \\
& & & & -\tilde K/2 & i\omega_1 &  \tilde K/2 \\
& & & & & -\tilde K/2 & i \omega_1
\end{bmatrix}
 \begin{bmatrix}
 w_{-\tilde N, N}\\
 w_{-\tilde N+2, N}\\
  \vdots\\
  \vdots\\
 w_{\tilde N-2, N}\\
 w_{\tilde N, N}
\end{bmatrix}
 \end{eqnarray*}
 with the state matrix being the product of
 a diagonal matrix and a  tridiagonal Toeplitz matrix,
 where $\tilde N=N$ for odd $N$ and $\tilde N=N-1$ for even order $N$.
Shown on the bottom plots  of Figure \ref{fig:kuramoto} are the approximation error
\begin{equation}\label{ENt.def}
 E(\theta_1(0), \omega_1, N, t)     =   \sup_{0\le s\le t}   \log_{10} \min\Big\{ 10,
  \max \big\{|w_{1, N}(s)- e^{i\theta_1(s)}|, 10^{-5}\big\}  \Big\}
 \end{equation}
of the finite-section approach \eqref{conventional.remark.eq1} to the  linearization \eqref{CarlemanFourier.eq0}  of the Kuramoto model  \eqref{kuramoto.def.one+} in the logarithmic scale, where
$\theta_1$ is the state in the dynamical system
 \eqref{kuramoto.def.one+}.
This suggests that $w_{1,N}$ may exhibit exponential convergence toward $e^{i\theta_1}$ over some time range,
and it has comparable approximation performance to the first block $v_{1, N}^\pm$ in  the finite-section approach to the proposed Carleman-Fourier linearization.
However, in contrast to the proposed Carleman-Fourier linearization framework in \eqref{CarlemanFourier.def}, we lack a rigorous mathematical proof for the exponential
convergence behavior observed in the above numerical simulation.

\section{Proof of Theorem \ref{main.thm}} \label{proof.section}

To establish the theorem, we begin by defining the error function $\pmb \eta_{k,N}(t) = {\bf y}_{k,N}(t) - {\bf y}_k(t)$ for $1 \leq k \leq N$. By leveraging \eqref{CarlemanFourier.def} and \eqref{Carleman.eq8}, the following relation holds:
\[
\dot{\pmb \eta}_{k,N}(t) - {\bf B}_{k,k} \pmb \eta_{k,N}(t) = \sum_{l=k+1}^N {\bf B}_{k,l} \pmb \eta_{l,N}(t) - \sum_{l=N+1}^\infty {\bf B}_{k,l} {\bf y}_l(t),
\]
with the initial condition $\pmb \eta_{k,N}(0) = {\bf 0}$ for $1 \leq k \leq N$.  Integrating this differential equation, we find:
\begin{eqnarray*}
\hspace{-0cm} e^{-{\bf B}_{k, k}t} \ {\pmb \eta}_{k,N}(t)
 &   = &     \int_0^t e^{-{\bf B}_{k, k} s} \Bigg( \sum_{l=k+1}^N {\bf B}_{k,l}{\pmb \eta}_{j, N}(s)
 -\sum_{l=N+1}^\infty {\bf B}_{k,l} {\bf y}_l(s) \Bigg) ds.
 \end{eqnarray*}
This representation highlights the role of higher-order terms and the iterative structure of the integral equation demonstrates how contributions from $\pmb \eta_{l,N}(t)$ for $l > k$ and the truncation at $N+1$ influence the dynamics of $\pmb \eta_{k,N}(t)$. Taking the $\ell^\infty$ norm of both sides, we have:
\begin{eqnarray}\label{etakn.infty}
   \|{\pmb \eta}_{k,N}(t)\|_\infty
& = & 
 \|e^{-{\bf B}_{k,k} t} \ {\pmb \eta}_{k,N}(t)\|_\infty\nonumber\\
 &    \le &
      \int_0^t
 \Big(\sum_{l=k+1}^N
 \|{\bf B}_{k,l} {\pmb \eta}_{l,N}(s)\|_\infty
  + \sum_{l=N+1}^\infty \|{\bf B}_{k,l} {\bf y}_l(s)\|_\infty\Big) ds\nonumber\\
 & \le  & 2
 Dk \int_0^t\sum_{l=k+1}^N r^{l-k} \|{\pmb \eta}_{l,N}(s)\|_\infty ds+
 \frac{2Dkr}{1-r}  r^{N-k} t,
\end{eqnarray}
where the first equality holds because ${\bf B}_{k,k}$ is a diagonal matrix with diagonal entries that are purely imaginary, and the second inequality follows from the observation that $\|{\bf y}_l(s)\|_\infty = 1$ and the fact that each row of ${\bf B}_{k,l}$ for $l \geq k+1$ contains at most two nonzero entries, each bounded by $Dk r^{l-k}$ as given by \eqref{fouriercoefficient.decay}.

 We define
\begin{equation}\label{ukn.def0}
u_{k,N}(t) = (1-r) r^{k-N-1} \|{\pmb \eta}_{k,N}(t/(2D))\|_\infty
\end{equation}
for $1 \leq k \leq N$. From this definition, it follows that
\begin{equation}\label{ukn.def}
u_{k,N}(t) \leq k \int_0^t \Big(1 + \sum_{l=k+1}^N u_{l,N}(s)\Big) ds
\end{equation}
for $1 \leq k \leq N$. We now establish the following critical estimate:
\begin{equation}\label{ukn.estimate}
1 + \sum_{l=k}^N u_{l,N}(t) \leq \sum_{m=0}^{N-k+1} {\binom{N}{m}} {\binom{N-k+1}{m}} t^m \quad {\rm for} \ t \geq 0,
\end{equation}
using induction on $k = N, \ldots, 2$.  Taking $k = N$ in \eqref{ukn.def}, we have
\begin{equation} \label{ukN.estimate}
1 + u_{N,N}(t) \leq 1 + N t,\end{equation}
 which confirms that \eqref{ukn.estimate} holds for $k = N$.  Next, we proceed by induction. Assume that the desired conclusion holds for $k+1$. Substituting this assumption into \eqref{ukn.def}, we obtain the result for $k$:
\begin{eqnarray*}   1+ \sum_{l=k}^N  u_{l, N}(t)
&
   \le  & \sum_{m=0}^{N-k} {\binom{N}{m}} {\binom{N-k}{m}}\big(t^m+ \frac{k}{m+1} t^{m+1}\Big)\\
    & = &  1+ \binom{N}{N-k} \frac{k}{N-k+1} t^{N-k+1}
    \\
     &   &
      + \sum_{m=1}^{N-k}
\Bigg( {\binom{N}{m}} {\binom {N-k}{m}}+ \frac{k}{m}
{\binom{N}{m-1}}{\binom{N-k}{m-1}}\Bigg)t^m\\
&  \le & \   1+ \binom{N}{N-k+1}  t^{N-k+1}
  + \sum_{m=1}^{N-k}
\Bigg( {\binom{N}{m}} {\binom {N-k}{m}}
  \nonumber\\
 &   &
\quad + \frac{N-m+1}{m}
{\binom{N}{m-1}}{\binom{N-k}{m-1}}\Bigg)t^m\\
& =&  \sum_{m=0}^{N-k+1}\binom{N}{m} \binom{N-k+1}{m} t^m.
    \end{eqnarray*}
Thus, the inductive proof can proceed as desired. Combining \eqref{ukn.estimate} for $k = 2$ with \eqref{ukn.def} for $k = 1$, we obtain
    \begin{eqnarray}
    u_{1, N}(t) & \le &    \sum_{m=0}^{N-1} \binom{N}{m} \binom{N-1}{m} \ \frac{t^{m+1}}{m+1} \nonumber
 =
   N^{-1} \sum_{m=1}^N \binom{N}{m-1}\binom{N}{m} \ t^m \nonumber\\
    & \hskip-0.08in \le & \hskip-0.08in
     N^{-1} t^{1/2} \sum_{m=1}^N \binom{N}{m-1}  t^{(m-1)/2} \Big(\sum_{n=1}^N \binom{N}{n} \ t^{n/2}\Big)\nonumber \\
   &  =  &
     N^{-1} t^{1/2}
    \left(\sum_{m=1}^N \binom{N}{m-1} t^{(m-1)/2}\right) \left(\sum_{n=1}^N \binom{N}{n} t^{n/2}\right)
    \le   N^{-1} t^{1/2} \left(1+t^{1/2}\right)^{2N}.
    \end{eqnarray}
This, together with \eqref{ukn.def0}, establishes \eqref{mainestimate} and thereby completes the proof of Theorem \ref{main.thm}.

\section{Conclusions and Discussions}

We introduced Carleman-Fourier linearization, a novel framework for analyzing nonlinear dynamical systems with quasi-periodic vector fields driven by multiple fundamental frequencies. By leveraging Fourier basis functions, this method transforms such systems into equivalent infinite-dimensional linear representations, offering a powerful new tool for their analysis. The efficacy of this approach was demonstrated through its application to the Kuramoto model, a widely studied system of coupled oscillators, where it achieved superior accuracy compared with the classical Carleman linearization technique.
 A key feature of our framework is the finite-section approximation, which includes quantifiable error bounds. We established that the primary block of the approximate solution converges exponentially to the solution of the original nonlinear system as the model size increases. These findings represent a significant advancement in the identification and analysis of quasi-periodic nonlinear systems, paving the way for efficient computational algorithms, rigorous reachability analysis, and applications across engineering, physics, biology, and quantum computing.

\bigskip
{\bf Acknowledgement}:\  The authors thank the anonymous reviewer for the insightful and constructive comments.

\end{document}